\documentclass[10pt, reqno]{amsart}
\usepackage{amssymb,latexsym}
\usepackage{eucal}
\newtheorem*{proposition*}{}
\newtheorem*{LemA}{Lemma A}

\newtheorem*{Lem10.4}{Lemma 10.4}
\newtheorem*{Lem17.3}{Lemma 17.3}
\newtheorem*{Lem18.3}{Lemma 18.3}
\newtheorem*{Lem18.5}{Lemma 18.5}
\newtheorem*{Lem19.0}{Lemma 19.0}
\newtheorem*{Lem19.1}{Lemma 19.1}
\newtheorem*{Lem19.1.3}{Lemma 19.1.3}

\newtheorem*{Th}{Theorem}
\newtheorem*{Cr1}{Corollary 1}
\newtheorem*{Cr2}{Corollary 2}

\newtheorem*{Proposition}{Proposition}

\DeclareMathOperator{\Hol}{\textup{Hol}}

\newcommand{\e}{{\varepsilon}}
\newcommand{\la}{{\langle}}
\newcommand{\ra}{{\rangle}}

\newcommand{\Y}{\mathcal{Y}}
\newcommand{\K}{\mathcal{K}}
\newcommand{\G}{\mathcal{G}}
\newcommand{\D}{\Delta}

\newcommand{\F}{\mathcal{F}}

\newcommand{\p}{\partial}
\newcommand{\E}{{\textup{E}}}
\newcommand{\kp}{\tau}

\begin{document}

\title[On HNN-extensions in the class of groups of large odd exponent]
{On HNN-extensions in the class of groups of large odd exponent}
\author{S. V. Ivanov}
\address{Department of Mathematics\\
University of Illinois \\
Urbana,   IL 61801, U.S.A. } \email{ivanov@math.uiuc.edu}
\thanks{Supported in part by NSF grant   DMS 00-99612}

\begin{abstract}
A sufficient condition for the existence  of HNN-extensions in the class of
groups of odd exponent $n \gg 1$ is given in the following form. Let $Q$ be a
group of  odd exponent $n > 2^{48}$ and  $\mathcal G$ be an HNN-extension of
$Q$. If $A \in \mathcal G$  then let $\mathcal F(A)$ denote the maximal
subgroup of $Q$ which is normalized by $A$. By  $\tau_A$ denote the
automorphism of $\mathcal F(A)$ which is induced by conjugation by $A$.
Suppose that for every $A \in \mathcal G$, which is not conjugate to an
element of $Q$, the group $\langle \tau_A, \mathcal F(A) \rangle$ has
exponent $n$ and, in addition, equalities $A^{-k} q_0 A^{k} = q_k$, where
$q_k \in Q$ and $k =0, 1, \dots, [2^{-16}n]$ ($[2^{-16}n]$ is the integer
part of $2^{-16}n$), imply that $q_0 \in \mathcal F(A)$. Then the group $Q$
naturally embeds in the quotient $\mathcal G / \mathcal G^n$, that is, there
exists an analog of the HNN-extension $\mathcal G$ of $Q$ in the class of
groups of exponent $n$.
\end{abstract}

\maketitle

\section{Introduction}

In this paper, we will generalize one of technical ideas of article
\cite{I02} which seems to be of independent interest and could be used for
future references.

Consider the following construction. Let $Q$ be a group of {\em exponent}
$n$ (that is, elements of $Q$ satisfy the identity $x^n \equiv 1$), $ \Y =
\{ y_1, \dots , y_m\} $ be an alphabet and $m \ge 1$.  Let  $P_{k,1}$,
$P_{k,2}$ be two  subgroups of $Q$, $k =1, \dots, m$, such that $P_{k,1}$ is
isomorphic to $P_{k,2}$ and
$$
\rho_k : P_{k,1} \to P_{k,2}
$$
be a fixed isomorphism. Then the standard HNN-extension $\G$ of $Q$ with
stable letters $y_1, \dots , y_m$ and isomorphisms $\rho_1, \dots , \rho_m$
is defined by the following relative presentation
\begin{equation}\label{1}
\G = \langle \   Q,  y_1, \dots , y_m    \  \|\
y_k p y_k^{-1} = \rho_k(p), \ p \in P_{k,1}, \   k =1, \dots, m \   \rangle .
\end{equation}

One might inquire whether there is a group $\K$ of exponent $n$ which
contains an isomorphic copy of $Q$  and some elements $ y_1, \dots , y_m$ so
that $y_k p y_k^{-1} = \rho_k(p)$ for all $p \in P_{k,1}$, $k =1, \dots, m$,
that is, there is an analog $\K$ of HNN-extension  of $Q$  with stable
letters $y_1, \dots , y_m$ and isomorphisms $\rho_k : P_{k,1} \to P_{k,2}$, $
k =1, \dots, m$, in the class of groups of exponent $n$. Clearly, the
existence of such a group $\K$ is equivalent to the natural embedding of $Q$
into the quotient  $\G / \G^n$.

It is also clear that in general the quotient $\G / \G^n$ need not contain
the natural copy of $Q$. For example, let $n$ be prime, $m=1$, let $P_{1,1}
= P_{1,2}$ be a subgroup  of order $n$  and $\rho_1(p) \neq p$, where $p \in
P_{1,1}$, $p \neq 1$. Then $y_1^n p y^{-n}_1 = \rho_1(p)  \neq p$,  whence
$p \in \G^n$. Here is another example. Let $P_{1,1} = P_{1,2} = Q$, that is,
$\rho_1 \in \mbox{Aut}Q$ and $\rho_1^n \neq 1$ in $\mbox{Aut}Q$. Then it is
clear that $Q \cap \G^n \neq \{ 1\}$.

The aim of this paper is to give a sufficient condition for the embedding $Q
\to \G / \G^n$. To do this, for every element $A \in \G$ we consider the
maximal subgroup $\F(A) \subseteq Q$ which is normalized by $A$ (clearly,
such an $\F(A) \subseteq Q$ is unique). For example, if $A \in Q$ then
$\F(A) = Q$. Let $\kp_A$ denote the automorphism of $\F(A)$ which is induced
by conjugation by $A$. Using this notation, it is easy to state and prove
(see Sect. 3) the following necessary condition for the embedding $Q \to \G
/ \G^n$.

\begin{Proposition}
Suppose that $Q$ is a group of exponent $n$, an HNN-extension $\G$ of $Q$ is
given by presentation $(1)$ and $Q$ naturally embeds in the quotient $\G /
\G^n$. Then  for every  $A \in \G$ the subgroup $\langle \kp_A, \F(A)
\rangle$ of $\Hol\F(A)$  has exponent $n$.
\end{Proposition}

A sufficient condition for the embedding $Q \to \G / \G^n$ given
in the following Theorem is our main result and is as close to the
necessary condition stated above  as we can get (in the statement
of Theorem, $[r]$ denotes the integer part of a real number~$r$).

\begin{Th}
Suppose that $Q$ is a group of odd exponent $n > 2^{48}$ and an
HNN-extension $\G$ of $Q$ is given by presentation $(1)$. Furthermore,
assume that for every $A \in \G$, which is not conjugate to an element of
$Q$, the subgroup $\langle \kp_A, \F(A) \rangle$ of $\Hol\F(A)$ has exponent
$n$ and, in addition, equalities $A^{-k} q_0 A^{k} = q_k$, where $q_k \in Q$
and $k =0, 1, \dots, [2^{-16}n]$, imply that $q_0 \in \F(A)$. Then the group
$Q$ naturally embeds in the quotient $\G / \G^n$.
\end{Th}

To prove this Theorem, we will make use of  the machinery of article
\cite{I94} (but first we will have to adjust it to the new situation which
is similar to  \cite{I02}). Interestingly, most of the machinery of article
\cite{I94}, created to solve the Burnside problem for even exponents $n \gg
1$, is "recycled" in this paper which, in particular, explains why we use
the estimate $n \ge 2^{48}$ of \cite{I94}.

As examples of application of this Theorem, we will state a couple
of immediate corollaries (other applications will be given
elsewhere). Note that these examples are similar to classical
applications of standard HNN-extensions proposed by  Higman,
B.~Neumann and H.~Neumann \cite{HNN49}.

\begin{Cr1}
Suppose that $G$ is a group of  odd exponent $n > 2^{48}$.
Then $G$ embeds into a group of exponent $n$  in which
every maximal cyclic subgroup has order $n$.
\end{Cr1}

\begin{Cr2}
Suppose that $G$ is a group of  prime exponent $n > 2^{48}$. Then $G$ embeds
into a group of exponent $n$ with $n$ classes of conjugate elements.
\end{Cr2}

Recall that a subgroup $H$ of a group $G$ is called {\em antinormal} if for
every $g \in G$ the inequality $gHg^{-1} \cap H \neq \{ 1\}$ implies that $g
\in H$. Note that earlier Mikhajlovskii \cite{M94} established the embedding
$Q \to \G / \G^n$, where odd $n \gg 1$, under the following assumptions: $m
=1$, the subgroups $P_{1, 1}$, $P_{1, 2}$ are both antinormal in $Q$ and $q
P_{1, 1} q^{-1} \cap P_{1, 2} = \{ 1\}$ for every $q \in Q$. We also remark
that this result follows from our Theorem (because, under these assumptions,
the subgroup  $\F(A)$ is trivial for every $A \in \G$ which is not conjugate
to an element of $Q$).

\section{Proving  Theorem}

According to the presentation (1), our basic alphabet is
$$
\Y_Q  = \{ y_1, \dots,   y_m, Q   \}
$$
and,  from now on  (unless stated otherwise), all words will be those in the alphabet
$
\Y_Q^{ \pm  1 } = \{  y_1, y_1^{-1}, \dots,   y_m, y_m^{-1} , Q   \} ,
$
called $\Y_Q$-{\em words}.
Let
$$
U_1  = S_{0,1} y_{k_{1,1}}^{\e_{1,1}} S_{1,1} \dots
y_{k_{\ell_1,1}}^{\e_{\ell_1,1}} S_{\ell_1,1} , \quad
U_2  = S_{0,2} y_{k_{1,2}}^{\e_{1,2}} S_{1,2} \dots
y_{k_{\ell_2,2}}^{\e_{\ell_2, 2}} S_{\ell_2,2} ,
$$
be  $\Y_Q$-words,  where
\begin{gather*}
\e_{1,1},  \dots, \e_{\ell_1, 1},  \e_{1,2}, \dots, \e_{\ell_2, 2} \in \{\pm 1 \} , \\
k_{1,1},  \dots, k_{\ell_1, 1},  k_{1,2}, \dots, k_{\ell_2, 2} \in \{1,\dots, m \} ,
\end{gather*}
$S_{0,1}, \dots,  S_{\ell_1,1}$ are $Q$-{\em syllables} of the word $U_1$
(that is, maximal subwords of $U_1$ all of whose letters are in $Q$; if
$S_{0,1}$ is, in fact, missing in $U_1$,  then we assume that $S_{0,1}=1 \in Q$ and
set $S_{\ell_1,1}=1$ if $S_{\ell_1,1}$ is not present in $U_1$)  and
$S_{0,2}, \dots, S_{\ell_2,2}$ are $Q$-syllables  of the word $U_2$.

We will write $U_1 = U_2$
if $\ell_1 = \ell_2$, $\e_{j, 1}= \e_{j, 2}$,
$k_{j, 1}= k_{j, 2}$ for all
$j =1, \dots, \ell_1$, and $S_{j, 1} = S_{j, 2}$ in the group $Q$
for all $j =0, \dots, \ell_1$.

The {\em length}  $|U_1|$ of a word $U_1$ is $\ell_1$, that is,
the number of occurrences of letters $y^{\pm 1}$ in $U_1$, where
$y \in \Y = \{ y_1, \dots, y_m \}$.

If  the (images of)  words $U_1, U_2$ are equal in the group
$\mathcal G(0) = \mathcal G$ given by the presentation
(1)  then we will write $U_1 \overset 0  = U_2$.

By induction on $i  \ge 0$  we will construct groups $\mathcal G(i)$. Assume that the
group $\mathcal G(i), i \ge 0$,  is already constructed as a quotient of
the group of  $\mathcal G(0)$
by means of defining relations. Define $\mathcal X_{i+1}$ to be  a maximal set of all
$\Y_Q$-words of length $i+1$ (if any) with respect to the following three properties.

\begin{enumerate}
\item[(AC1)] Every word $A \in \mathcal X_{i+1}$ begins with $y$ or $y^{-1}$, where
$y \in \Y = \{ y_1, \dots, y_m \}$.
\item[(AC2)] The image of every word
$A \in \mathcal X_{i+1}$ has infinite order in the group $\mathcal G(i)$.
\item[(AC3)]  If $A, B$ are distinct elements of  $\mathcal X_{i+1}$  then the
image of  $A^{n}$ is not conjugate in $\mathcal G(i)$ to the image of $B^{n}$ or
$B^{-n}$.
\end{enumerate}

Note  that it follows from the analog of Lemma 18.2 in rank $i \geq 0$ that
the set $\mathcal X_{i+1}$ is nonempty.

Similar to  \cite{O82a}, \cite{O89}, \cite{I94}, we will call a word
$A \in \mathcal X_{i+1}$ a  {\em period of rank} $i+1$.

Now we define the group $\mathcal G(i+1)$  by imposing all relations
$A^{n} =1$, $A \in  \mathcal X_{i+1}$,
called {\em  relations of rank $i+1$}, on the group  $\mathcal G(i)$:
$$
\mathcal G(i+1) =
\langle \; \mathcal G(i) \; \| \;  A^{n} =1, \ A \in \mathcal X_{i+1} \;  \rangle .
$$

It is clear that
$$
\mathcal G(i+1) = \langle \; \mathcal G  \; \| \;  A^{n} =1, \
A \in \cup_{j=1}^{i+1} \mathcal X_{j} \;  \rangle .
$$

We also define the limit group $\mathcal G(\infty)$  by imposing on the free product
$$
Q * \langle y_1 \rangle_\infty *
\dots  *  \langle y_m  \rangle_\infty
$$
of all relations of all ranks $i = 0, 1, 2, \dots$:
\begin{multline} \label{2??}
\mathcal G(\infty) = \langle \  Q, \; y_1, \dots , y_m   \   \|   \
y_k p  y_k^{-1}= \rho_k(p) ,  \; p \in P_{1,k}, \ k =  1, \dots, m , \\
A^{n} =1, \  A \in \cup_{j=1}^{\infty} \mathcal X_{j} \  \rangle .
\end{multline}

The main technical result relating to the group
$\mathcal G(\infty)$ is the following.
(Observe that Lemma A obviously implies Theorem.)

\begin{LemA} Suppose that
the hypothesis of Theorem holds.  Then the
group $\mathcal G(\infty)$ given by  presentation $(2)$
is naturally isomorphic to $\G / \G^n$ and the group
$Q$ naturally embeds in $\mathcal G(\infty)$.
\end{LemA}

We will make use of the machinery of article \cite{I94} to prove
Lemma A (more applications of this machinery can be found in
\cite{IO96}, \cite{IO97}, \cite{I00}, \cite{I02}).
In order to do this  we will have to make necessary changes in definitions,
statements of lemmas of \cite{I94}  and their proofs.

First of all, diagrams over the group $\mathcal G(0) = \mathcal G$ given by
presentation (1)  (or, briefly, over $\mathcal G$),
 called  {\em diagrams of rank} 0, are defined to be maps
that have two types of  2-cells.

A 2-cell $\Pi$ of the first type, called a 0-{\em square}, has four edges in
its counterclockwise oriented boundary (called {\em contour})
$
\p\Pi = e_1 e_2 e_3 e_4
$
and
$$
\varphi(e_1) = \varphi(e_3)^{-1} = y_k , \ \   \varphi(e_2) = g,  \ \
\varphi(e_4) = {\rho_k}(g)^{-1} ,
$$
where, as in  \cite{I94},   $\varphi$ is the labelling function, $g \in P_{k,1}$
(perhaps, $g=1$)  and $y_k \in \{y_1, \dots, y_m \}$. Observe that we use Greek letters
with no indices exactly as in \cite{I94}  (in particular, see table (2.4) in \cite{I94}).

A 2-cell $\Pi'$ of the second type, called a 0-{\em circle},  has $\ell \ge 2$ edges
in its contour  $\p\Pi' = e_1 \dots e_\ell$ so that
$\varphi(e_1), \dots ,  \varphi(e_\ell) \in Q$ and the word
$$
\varphi(\p\Pi' ) = \varphi(e_1)  \dots  \varphi(e_\ell)
$$
equals 1 in $Q$.

Note that this definition of  2-cells in a diagram of rank 0 is analogous
to the corresponding definitions of  \cite{I02}.

Furthermore, we take into account that there are multiple periods of rank $i$ and the symbol
$A_i$ will denote one of many periods of rank $i$ (note that the length $|A_i|$
of  $A_i$ is now $i$).

In the definition of an $A$-periodic word, it is now assumed that $A$ starts with
$y$ or $y^{-1}$, where $y \in \Y$,
and $A$ is not conjugate in $\G$  to a power $B^\ell$ with $|B| < |A|$.

In addition to cells of positive rank,
we also have (as in \cite{O89}, \cite{IO96}, \cite{I02}) cells of rank 0
(which are now either 0-squares or 0-circles).
The equality $r(\D)= 0$ now means that all cells in $\D$ have rank 0.

If $ p = e_1 \dots e_\ell , $ where $e_1, \dots, e_\ell$ are edges, is a path in a
diagram $\D$ of rank $i$ (that is, a diagram over the group  $\mathcal G(i)$) then the
$y$-length $|p|$ of $p$ is $|\varphi(p)|$, that is, $|p|$ is the number of all edges of
$p$ labelled by $y^{\pm 1}$, where $y \in \Y$. The (strict) length of  $p$, that is,  the
total number of edges of $p$, is $\ell$ and denoted by $\| p \|$.

In the definition (A1)--(A2)  of $j$-compatibility (p.13 \cite{I94}) we eliminate
the part (A2)  because $n$ is odd and,
similar to \cite{I94} in the case when $n$ is odd,
it will be proven in a new version of Sect. 19 \cite{I94}  that
there are no $\mathcal F(A_j)$-involutions and there is no  $j$-compatibility of
type (A2).

We can also drop the definition of self-compatible cells (p.13 \cite{I94}) because they
do not exist when $n$ is odd (which is again analogous to \cite{I94} in the case when $n$
is odd). Thus all lemmas in \cite{I94} whose conclusions deal with self-compatible cells,
compatibility of  type (A2) actually claim that their assumptions are false (e.g.,  see
lemmas of Sect. 12 \cite{I94}). On the other hand, the existence of self-compatible cells
in assumptions of lemmas of \cite{I94}  is now understood as the existence of
noncontractible  $y$-annuli  which we are about to define.

A $y$-{\em annulus}, where $y \in \Y$, is defined to be an annular
subdiagram $\Gamma$ in a diagram $\D$ of rank $i$
such that $\Gamma$ consists of 0-squares $S_1, \dots, S_k$ so that if
$$
\p S_\ell = f_{1, \ell} e_{1, \ell} f_{2, \ell} e_{2, \ell} ,
$$
where $e_{1, \ell}, e_{2, \ell}$ are $y$-{\em edges} (that is, labelled by $y$ or
$y^{-1}$, where $y \in \Y$) of $\p S_\ell$, $1 \le \ell \le k$, then  $e_{2, \ell}= e_{1,
\ell+1}^{-1}$, where the second subscript is  $\mbox{mod} k$, $\ell =1, \dots, k$. If
$\Gamma$ is contractible into a point  in $\D$ then we will call $\Gamma$ a {\em
contractible} $y$-annulus. Otherwise,  $\Gamma$ is a  {\em noncontractible} $y$-annulus.
If   $\Gamma$ is contractible in $\D$ and bounds a simply connected subdiagram $\Gamma_0$
in $\D$ with $\varphi( \p \Gamma_0) = 1$ in $Q$  then $\Gamma$ is termed a  {\em
reducible} $y$-annulus.

In the definition of a reduced (simply connected or not)  diagram $\D$
of rank $i$ (p.13 \cite{I94}),  we additionally require that $\D$ contain no
reducible $y$-annuli.

As in \cite{I94}, we can always remove reducible pairs and reducible $y$-annuli
in a diagram $\D$ of rank $i$ to obtain from $\D$ a reduced
diagram $\D'$ of rank $i$. Note that in general  it is not possible to get rid of noncontractible
$y$-annuli  (in non-simply connected  diagrams of rank $i$).

In the definition of a 0-bond $\E$  between $p$ and $q$ (p.15 \cite{I94}) we now require that
$\E$ consist of several  0-squares $S_1, \dots, S_k$ so that if
$$
\p S_\ell = f_{1, \ell} e_{1, \ell} f_{2, \ell} e_{2, \ell} ,
$$
where $e_{1, \ell}, e_{2, \ell}$ are $y$-edges  of $\p S_\ell$,
$1 \le \ell \le k$, then $e_{1,1}^{-1} \in p$, $e_{2,\ell}= e_{1,\ell+1}^{-1}$,
$\ell =1, \dots, k-1$, and $e_{2, k}^{-1} \in q$.

The {\em standard} contour of the 0-bond $\E$ between $p$ and $q$  is
$$
\p \E = e_{1, 1}^{-1}(f_{1, 1}^{-1} \dots f_{1, k}^{-1}) e_{2, k}^{-1}
(f_{2, k}^{-1} \dots f_{2, 1}^{-1})
$$
and the edges
$e_{1,1}^{-1}, e_{2,k}^{-1}$ are denoted by $\E \wedge p$, $\E \wedge q$,
respectively.

In the definition of  a simple in rank $i$ word $A$ (p.19  \cite{I94}),
we additionally require that $|A| >0$ and  $A$ start with $y^{\pm 1}$, where $y \in \Y$.
Observe that it follows from Lemma 18.2 (in rank $i-1$)  and definitions that
a period of rank $i$ is simple in rank $i-1$ (and hence in any rank $j \leq i-1$).

In the definition of a tame diagram of rank $i$ (p.19 \cite{I94}), we make two changes.
First, in  property (D2), we require that if 0-squares $S_1, \dots, S_k$ form
a subdiagram $\E$ as in  the definition of a 0-bond and $p= q = \p\Pi$,
where $\Pi$ is a cell of rank $j >0$ in $\D$,  then
$\E$ is a 0-bond between  $\p\Pi$ and $\p\Pi$ in $\D$. Second, we add the
following property.
\begin{enumerate}
\item[(D3)]  $\D$ contains no contractible $y$-annuli.
\end{enumerate}

In the definition of  a complete system (p.23 \cite{I94}) we require in (E3) that $e$ be a
$y$-edge.

In Lemma 4.2, the strict length $\| s_1 \|,  \| s_2 \|$ of $s_1, s_2$ is meant.

In the definition of  the weight function $\nu$ (p.28 \cite{I94}),
we require in (F1) that $e$ be a $y$-edge.
In (F2), we allow that $e$ is not  a $y$-edge.

In the beginning of the proof of Lemma 6.5, we note that the lemma is obvious if
$\D$ contains no cells of positive rank. In general, repeating  arguments of \cite{I94},
we always understand "cells" as cells of positive rank and keep in mind
the existence of cells of rank 0.

In the definition of the height $h(W)$ of a word $W$ (p.89 \cite{I94}),  we additionally set
$
h(W) = \tfrac{1}{2}
$
if  $W \overset i  \neq 1$  and
$W$ is conjugate in rank $i$ to a word $U$ with $|U| = 0$
(that is, $U \in Q$ and $U  \overset i  \neq 1$).

In Lemma 10.2, we allow the extra case when $h(W) = \tfrac{1}{2}$.

Here is a new version of Lemma 10.4.

\begin{Lem10.4}  $(a)$  If a word $W$ has finite order $d > 1$ in the group
$\mathcal G(i)$ then $n$ is divisible by $d$.

$(b)$ Every word $W$  with  $|W| \leq i$  has finite order in rank $i$.
\end{Lem10.4}

{\em Proof.} (a)  By  Lemma 10.2, either  $h(W) = \tfrac{1}{2}$
(and then our claim is immediate from  $Q$ being of exponent $n$)
or, otherwise, $W$ is conjugate in rank $i$
to a word of the form $A^k F$, where  $A$ is a period of rank $j \le i$,
$0 < k < n$  and $F \in \mathcal F(A)$.
In the latter case, it follows from Lemma 18.5(c) in rank $j-1 <i$ that
$
(A^k F)^{n} \overset {j-1}  = A^{k n} .
$
Therefore,
$
(A^k F)^{n} \overset {j-1}  = A^{k n}  \overset {j}  = 1 ,
$
whence $W^{n} \overset i = 1$ as required.

(b)  By induction, it suffices to show that  every word $W$ with
$|W| = i$ has finite order in rank $i$ (for $i=0$ this is obvious).
It follows from the definition of periods
of rank $i \geq 1$ that if $W$ has infinite order in rank $i-1$ then
$W^{n}$ is conjugate in rank $i-1$ to $A^{\pm n}$, where
$A$ is a period of rank $i$.
Therefore, $W^{ n} \overset i = 1$ as desired.

Lemma 10.4 is proved. \qed
\smallskip

In Lemma 10.8, we drop part (b)  of its conclusion (and keep in mind that the term
"reducible cell" now means "$y$-annulus"). Note that the height
of $\varphi(q_1t_1)$ in Lemma 10.8 is at least $1$ hence noncontractible $y$-annuli
in $\D_0$, $\D_0^r$ are impossible (for otherwise, the height  $h(\varphi(q_1t_1))$
of  $\varphi(q_1t_1)$ would be at most $\tfrac {1}{2}$).

Lemma 10.9 is no longer needed
for no path $q$ is (weakly) $j$-compatible with itself.

In the definition of a $U$-diagram of rank $i$ (p.134 \cite{I94}), we allow in property
(U3) that the height $h(\varphi(e))$  of  $\varphi(e)$  is  $\tfrac {1}{2}$.

Lemma 12.3 now claims that there are no  $U$-diagrams of rank $i$. Recall that this
agrees with our convention that if the conclusion  states the existence of self-compatible
cells or $j$-compatibility of type (A2) then the assumption is false.

The analogs of Lemmas 13.1--16.6 are not needed.

In the hypothesis of  Lemma 17.1,  we now suppose that
one can obtain  from $\D_0$ an annular reduced diagram
of rank $i$  which contains no noncontractible $y$-annuli  by
means of  removal of  reducible pairs and reducible  $y$-annuli.

According to our convention, in the statement of Lemma 17.2, we replace
the  phrase "one has to remove a reducible cell  to reduce $\D_0$" by
"one encounters a noncontractible $y$-annulus when reducing $\D_0$".
In the conclusion of Lemma 17.2 and in its proof, we disregard reducible cells,
$\mathcal F(A_j)$-involutions and consider, instead, noncontractible $y$-annuli
and their 0-squares.

The new version of Lemma 17.3 is stated as follows.

\begin{Lem17.3} Let $\D$ be a disk reduced diagram of rank $i$
whose contour is $\p\D =bpcq$, where $\varphi(p)$ and $\varphi(q)^{-1}$ are
$A$-periodic words and $A$ is a simple in rank $i$ word with $|A| = i+1$
(in particular, $A$ is a period of rank $i+1$).
Suppose also that $\D$ itself is a contiguity
subdiagram between $p$ and $q$ with $\min(|p|,|q|) > L|A|$.
Then, in the notation of Lemmas 9.1--9.4,   there
exists a rigid subdiagram of  the form $\D(m_1, m_2)$  in $\D$ $( = \D(1, k))$
such that
$$
r(\D(m_1, m_2))=0
$$
and the following analogs of  inequalities (17.25)--(17.26)  (p.222
\cite{I94}) hold.
\begin{align*}
|q(m_2,m_1)| & > |q(k,1)|-4.4|A| > (L-4.4)|A| , \\
|p(m_1,m_2)| & > |p(1,k)|-4.4|A| > (L-4.4)|A|  , \\
|x_{m_1}|  & = |y_{m_2}| =0  .
\end{align*}
\end{Lem17.3}

{\em Proof.} To prove this new version of Lemma 17.3,
we repeat the argument of the  beginning of the proof of Lemma 17.3 \cite{I94}. As there,
making use of Lemma 17.2, we prove Lemma 17.3.1. After that,
arguing as in the proof of  Lemma 17.3.2, it is easy to show, using Lemma 12.3,
that
$$
r(\D(m_1, m_2))=0 ,
$$
as required.  \qed
\smallskip

When proving the analog of  Lemma 18.2, we pick a word $B= B(a_1, a_2)$ in
the  alphabet  $\{ a_1, a_2\}$ of length $i+1$ so that $B$ has the same
properties as those  in \cite{I94} and, in addition, first and last letters
of $B$ are distinct (the existence of such a word easily follows from Lemma
1.7 \cite{I94}). Next, consider a word $B(a_1, a_2, q)$ in the alphabet $\{
a_1, a_2, q \}$ which is obtained from $B(a_1, a_2)$ by plugging in an
element $q \in Q, q \neq 1$, between each pair of consecutive letters of the
word $B(a_1, a_2)$. Then we replace each occurrence of the letter $a_1$ in
$B(a_1, a_2, q)$ by $y$, $y \in \Y$,  and  each occurrence of the letter
$a_2$ in $B(a_1, a_2, q)$ by $y^{-1}$. Clearly, we have a word $B(y, y^{-1},
q)$ with $|B(y, y^{-1}, q)|= i+1$. Now, in view of Lemmas 10.2, 10.4, we can
repeat the arguments of the proof of Lemma 18.2 without any changes.
\smallskip

Let $A$ be a period of rank $i+1$.   As in Sect. 1, by
$$
\mathcal F(A)
$$
denote a maximal subgroup of  $Q \subseteq \mathcal G(0)$
with respect to the property that $A$ in rank $0$
normalizes this subgroup $\mathcal F(A)$.
Observe that  $Q$ naturally embeds in $\mathcal G(i)$ by Lemma 6.2
and so we can also consider $\mathcal F(A)$ as a subgroup of
$\mathcal G(i)$.

Here is a new version of Lemma 18.3.

\begin{Lem18.3}  Suppose that $A$ is a period of rank
$i+1$. Furthermore, let $\D$ be a disk reduced diagram of rank $i$
such that  $\p\D = b p c q$, where $\varphi(p), \;
\varphi(q)^{-1}$ are $A$-periodic words with $\min(|p|,|q|) >
\tfrac 1 2 \beta n |A|$, and $\D$  itself  be a contiguity
subdiagram between  sections $p$ and $q$. Then there exists a
$0$-bond $\E$ in $\D$ with the standard contour $\p \E = b_{\E}
p_{\E} c_{\E} q_\E$, where $p_\E = \E \wedge p$, $q_\E = \E \wedge
q$, $(p_\E)_-$, $(q_\E)_+$ are  phase vertices  of  $p, q$,
respectively,  such that
$$
\varphi(b_\E) \in  \mathcal F(A) ,  \qquad
 A^{n} \varphi(b_\E)  A^{-n} \overset i  =  \varphi(b_\E)
$$
and for every  integer $k$ one has
$$
(A^k \varphi(b_\E))^{n}  \overset i  =  A^{k n } .
$$
\end{Lem18.3}

{\em Proof.} Lemma 17.3 enables us to assume that
$$
r(\D) = 0, \quad  \min( |p|, |q| ) > (  \tfrac 1 2 \beta n  - 5) |A| , \quad  |b| = |c|=0 .
$$
In particular, there are $| p|$ 0-bonds between $p$ and $q$ in $\D$.
Let $\E$ be a 0-bond between sections $p$ and $q$ and
$$
\p \E = b_\E p_\E c_\E q_\E
$$
be the standard contour of $\E$,
where $p_\E = \E \wedge p$, $q_\E = \E \wedge q$.
It is clear that $\mbox{div}((p_\E)_-, (q_\E)_+)$
does not depend on $\E$. Suppose that
\begin{equation} \label{?}
\mbox{div}((p_\E)_-, (q_\E)_+) \neq 0 .
\end{equation}

$A$-Periodically extending $p$ or $q^{-1}$ on the left as in the beginning of the  proof
of Lemma 17.3 (see Fig. 17.4(a)--(b) in \cite{I94}), we will get a diagram $\D'$ with
$\p \D' = b' p' c' q'$ such that both $\varphi(p')$ and  $\varphi(q')^{-1}$ begin with a
cyclic permutation $\bar A$ of $A$ such that $\bar A$ starts with
$y^{\pm 1}$, where $y \in \Y$.

As in the  proof of  Lemma 17.1,  we can easily get,  making use of (3),  that the
annular diagram $\D'_0$ (obtained from $\D'$ as in Lemma 17.1)
contains no $y$-annuli, in particular, $\D'_0$ is already reduced.
Therefore,  Lemma 17.1 applies to $\D'$ and yields that
$\varphi(b') \overset i  = 1$. It follows from (3) that $|b'| >0$ and so $\bar A$
is not cyclically reduced in rank $i$. This contradiction proves that (3) is false and
$
\mbox{div}((p_\E)_-, (q_\E)_+)  =  0 .
$

Without loss of generality, we may assume that words
$\varphi(p)$, $\varphi(q)^{-1}$ start with $A$ and  the word $A$
starts with $y^{-1}$, where $y \in \Y$.

Using the notation of Lemma 9.1, let $\E_1, \dots, E_{|p|}$ be all
(consecutive along $p$)  0-bonds between $p$ and $q$ with standard contours
$$
\p \E_\ell = b_\ell  p_\ell  c_\ell  q_\ell ,
$$
where $p_\ell  = \E_\ell \wedge p$, $q_\ell  = \E_\ell \wedge q$
and  $1 \le \ell \le |p|$.

Next,  consider  words
$$
V_0 = \varphi(b_1), \ V_1 = \varphi(b_{|A| +1}), \dots ,
V_t =   \varphi(b_{t|A| +1}),  \dots ,
V_{[ 2^{-16} n] +1}=   \varphi(b_{ ([ 2^{-16} n]+1)|A| + 1}) .
$$
Recall that $|p| > (\tfrac 1 2 \beta n - 5)|A| >  ([ 2^{-16} n] +2)|A|$.

It is clear that
$$
A^{-1} V_t A = V_{t+1}
$$
and $V_t, V_{t+1}$ are words in $Q$, where $t = 0, \dots , [ 2^{-16} n] $.

By definitions, $A$ is simple in rank $i$ and so $A$ is not conjugate in rank
$i$ to an element of $Q$.
Therefore, it follows from Theorem's hypothesis that the word $V_0 = \varphi(b_1)$
belongs to the subgroup $\F(A) \subseteq Q$. Since the group $\langle \kp_A, \F(A) \rangle$
has exponent $n$, it follows that
$$
A^n \phi(b_1) A^{-n}   \overset 0  =  \kp_A^n(  \varphi(b_1))  \overset 0  =    \varphi(b_1)
$$
and
\begin{multline*}
(A^k  \phi(b_1)   )^{n}   \overset 0  =  (A^k  V_0  )^{n}   \overset 0  =
A^{k} V_0 A^{-k} A^{2k} V_0 A^{-2k} \dots
A^{k n} V_0 A^{-k n}  A^{k n}  \overset 0  = \\
\overset 0  = \kp_A^{k}  V_0  \kp_A^{-k} \kp_A^{2k}  V_0  \kp_A^{-2k} \dots
 \kp_A^{n k}  V_0  \kp_A^{- n k}  A^{k n}  \overset 0  =
(\kp_A^{k}  V_0)^n   \kp_A^{- n k}  A^{k n}  \overset 0  = A^{k n} ,
\end{multline*}
as required.  Lemma 18.3 is proven. \qed
\smallskip

Lemma 18.4 is not needed.

\begin{Lem18.5}
Let $A$ be a period of rank $i+1$ and  $\mathcal F(A)$ be the  maximal subgroup of
$Q \subset \mathcal G(0)$ with respect to the property that $A$
normalizes $\mathcal F(A)$ in rank $0$.  Then the following are true.

$(a)$  The subgroup $\mathcal F(A)$  is defined uniquely.

$(b)$   Suppose $\D$ is a disk reduced diagram of rank $i$ with $\p\D= bpcq$,
where $\varphi(p)$, $\varphi(q)^{-1}$ are $A$-periodic words with
$\min(|p|,|q|) > \tfrac 1 2 \beta n |A|$, such that $\D$ itself is a contiguity subdiagram
between $p$ and $q$. Then there is a $0$-bond $\E$ between $p$ and $q$
with the standard contour $\p \E = b_\E p_\E c_\E  q_\E$,
where $p_\E = \E \wedge p$,
$q_\E = \E \wedge q$, such that $(p_\E)_-$,  $(q_\E)_+$ are phase vertices of $p, q$,
respectively, and $\varphi(b_\E) \in \mathcal F(A)$.

$(c)$    $A^{n}$ centralizes  the subgroup $\mathcal F(A)$ and
if $F \in \mathcal F(A)$,  $k$ is an  integer
then
$$
(A^k F)^{n} \overset i  = A^{k n} .
$$

$(d)$   The subgroup $\langle \mathcal F(A), \; A \rangle$ of $\mathcal G(i)$ has
the property that a word $X$ belongs to $\langle \mathcal F(A), \; A \rangle$
if and only if there is an integer $m \ne 0$ such that
$X A^m X^{-1} \overset i = A^m$.
\end{Lem18.5}

{\em Proof.} (a)  This is obvious from definitions.

(b)  This follows from Lemma 18.3.

(c) These claims can be proved as similar claims of Lemma 18.3.

(d) By part (c), it suffices to show that an equality
$$
X  A^{m} X^{-1} \overset i  = A^{m} ,
$$
where $m \neq 0$,  implies that $X \in  \langle A,   \mathcal F(A) \rangle$.
Arguing  exactly as in the proof of part (c) of Lemma 18.5 \cite{I02},
we can show that
$$
X \in \langle A \rangle \varphi(b_\E)  \langle A \rangle \subseteq \mathcal G(i) ,
$$
where $\varphi(b_\E) \in  \mathcal F(A)$ by Lemma 18.3. Thus
$X \in  \mathcal F(A)$  and Lemma 18.5 is proven.  \qed
\smallskip

Let us state a new version of Lemma 19.1.

\begin{Lem19.1}  There is no disk diagram  $\D$ of rank $i$ such that
$\p \D = b p c q$, where $p$, $q$ are
$A$-periodic sections with $|p|, |q|  > \tfrac 13 {\beta} n |A|$, $A$ is a period
of rank $i+1$, and $\D$ itself is a contiguity subdiagram between $p$ and $q$.
\end{Lem19.1}

{\em Proof.} Arguing on the contrary, we assume the existence of such a diagram $\D$
and,  replacing the coefficient $N$  ($N = 484$ as defined in  (17.1) on p.212  \cite{I94})
by $\tfrac 13 {\beta} n$, repeat the proof of
Lemma 19.1 \cite{I94}  up to getting equality (19.23) (p.290 \cite{I94}) which now reads
\begin{equation} \label{16}
\varphi(d) A^{n} \varphi(d)^{-1} \overset i  =  A^{-n} .
\end{equation}

\begin{Lem19.1.3}  The equality $(4)$ is impossible.
\end{Lem19.1.3}

{\em Proof.} Arguing on the contrary, we note that it follows from (4) that
\begin{equation} \label{17}
\varphi(d)^2 A^{n} \varphi(d)^{-2} \overset i  =  A^{n} .
\end{equation}

Recall that $|A| = i+1$ and, by Lemma 19.1.1, $\varphi(d)$ is conjugate  in rank $i$
to a word $W$ with
$$
|W| < (0.5+\xi )|A| <|A| .
$$
Hence, by Lemma 10.4(b),   $\varphi(d)$ has finite order  in rank $i$ and,
by Lemma 10.4(a),
$$
\varphi(d)^{n} \overset i = 1 .
$$
This means  that
$\varphi(d)  \in \langle \varphi(d)^2 \rangle$ in rank $i$ and so equalities
(4) and (5) imply that
$$
A^{2n}   \overset i  = 1 .
$$

A contradiction  to the definition of a period of rank $i+1$ completes  proofs
of  Lemmas 19.1.3 and 19.1. \qed
\smallskip

Analogues of Lemmas 19.2--19.6 are no longer needed.

The statements and proofs of  Lemmas 20.1--20.2 are retained.
\smallskip

Having made all necessary changes, we can now turn to the group $\mathcal G(\infty)$
given  by presentation (2).

It follows from Lemma 10.4(b) that every word $W$ has finite order in rank $i  \geq |W|$.
Then, by Lemma 10.4(a), $W^{n}  \overset i  = 1$ provided that $i \geq |W|$. Thus, the
group $\mathcal G(\infty)$ has exponent $n$. Now it is clear that the group $G(\infty)$
is naturally isomorphic the quotient $\G / \G^n$.

Suppose that $W$ is a word with $|W| = 0$, that is, $W \in Q$. Let $W = 1$ in
the group $\mathcal G(\infty)$. Then there is an $i$ such that $W \overset i  = 1$.
Consider  a reduced diagram $\D$ of rank $i$ for this equality. Since $|\p \D | =0$,
it follows from Lemma 6.2 that $r(\D) =0$, that is, $\varphi(\p \D) =1$ in $Q$.
Thus $Q$ naturally embeds in $\mathcal G(\infty)$.

The proofs of Lemma A and Theorem are  complete. \qed

\section{Proofs of Proposition  and Corollaries}

First we  prove Proposition.
Arguing on the contrary, assume that for some $A \in \G$
the subgroup $\langle\kp_A, A \rangle$
of $\Hol\F(A)$ has no exponent $n$. Then either $\kp_A^n \neq 1$
in $\mbox{Aut}\F(A)$ or  $\kp_A^n = 1$ and there is a $q \in Q$ so that
$( \kp_A^k q)^n \neq 1$ in  $\Hol\F(A)$ for some $k$.

In the first case,  there are  $q_1, q_2$ in $Q$ such that
$q_1 \neq q_2$ and $A^n q_1 A^{-n} = q_2$. Hence, in view of $A^n \in \G^n$,
we have that $q_1^{-1} q_2 \in \G^n$, contrary to $\G^n \cap Q = \{ 1\}$.

Now suppose that $\kp_A^n = 1$ and
$( \kp_A^k q)^n \neq 1$ in  $\Hol\F(A)$ for some $q \in Q$ and $k$.
Then
\begin{multline*}
(A^k q  )^{n}     =
A^{k} q A^{-k} A^{2k} q A^{-2k} \dots
A^{k n} q A^{-k n}  A^{k n}   = \\
= \kp_A^{k}  q  \kp_A^{-k} \kp_A^{2k}  q  \kp_A^{-2k} \dots
 \kp_A^{n k} q  \kp_A^{- n k}    A^{k n}  =
(\kp_A^{k}  q)^n   \kp_A^{-n k}  A^{k n} .
\end{multline*}
Then it follows from $ \kp_A^{-n k} = 1$
and $( \kp_A^k q)^n \neq 1$ that
$
(A^k q)^n = q' A^{nk} ,
$
where $q' \in Q$, $q' \neq 1$.  Hence, $q' \in \G^n$,
contrary to $\G^n \cap Q = \{ 1\}$. Proposition is proved. \qed
\smallskip

Let us prove Corollary 1. Let $G$ be a group of odd exponent $n > 2^{48}$. By the
standard induction argument (see similar applications of HNN-extensions in \cite{HNN49}
or \cite{LS77}), it suffices to show that if $g \in G$ then $G$ embeds in a group $H$
such that $H$ has exponent $n$ and $H$ contains an element $z$ of order $n$ with $g \in
\la z \ra$. To do this, consider a cyclic group $\langle x \rangle_n$ of order $n$
generated by $x$ and let
$$
Q = G *_n \langle x \rangle_n = G * \langle x \rangle_n /
( G * \langle x \rangle_n  )^n
$$
be the free Burnside  $n$-product of $G$ and $\langle x \rangle_n$ (which is
the quotient of the free product $G * \langle x \rangle_n$ by $( G * \langle
x \rangle_n  )^n$). Recall that, by Lemma 34.10 \cite{O89}, if $n$ is odd and
$n \gg 1$ (say, $n
> 10^{10}$) then the factors of a free Burnside $n$-product are antinormal in the
product.

Let $x^k$ have the same order as that of $g$ and consider the following
HNN-extension of $Q$
$$
\G = \langle \;  Q, y \; \|  \; y g y^{-1} = x^k \;  \rangle .
$$
Assume that $A \in \G$, $A$ is not conjugate in $\G$ to an element of $Q$, and $A^{-1}
q_0 A = q_1$, $A^{-1} q_1 A = q_2$ in $\G$, where $q_0, q_1, q_2 \in Q$. Then it is not
difficult to see that for every $\ell >0$ it is true that $A^{-\ell} q_0 A^{\ell} \in Q$
and, more specifically, $A^{-\ell} q_0 A^{\ell} = q_A^{-\ell} q_0 q_A^{\ell}$, where $q_A
\in Q$ is the result of deletion in $A$ of all occurrences of $y^{\pm 1}$ and
$x$-syllables (we assume here that $A$ is a reduced in $\G$ word, that is, $|A|$ is
minimal).  Therefore, $q_0 \in \F(A)$, $\kp_A^n =1$ and
\begin{multline*}
( \kp_A^{-\ell} q_0 )^{n}     = \kp_A^{-\ell} q_0 \kp_A^{\ell} \kp_A^{-2\ell} q_0
\kp_A^{2\ell} \dots \kp_A^{-n\ell} q_0 \kp_A^{n\ell}  \kp_A^{-n\ell} = \\
= q_A^{-\ell} q_0 q_A^{\ell} q_A^{-2\ell} q_0 q_A^{2\ell} \dots q_A^{-n\ell} q_0
q_A^{n\ell}  \kp_A^{-n\ell} = (q_A^{-\ell} q_0)^n q_A^{n \ell} \kp_A^{-n\ell} = 1 .
\end{multline*}
Since $q_0$ is an arbitrary element of $\F(A)$, this implies that the group $\la \kp_A,
\F(A) \ra$ has exponent $n$ and so Theorem applies to $\G$. Hence, $Q$ naturally embeds
in $\G / \G^n$ and $(y^{-1} x y)^k = g$ in $\G / \G^n$, as required. Corollary 1 is
proved. \qed
\smallskip

The proof of  Corollary 2 is similar. Let $Q$ be a group of prime exponent $n
> 2^{48}$. As before, by induction,  it suffices to show that if $\langle p_1
\rangle$, $\langle p_2 \rangle$ are two nontrivial cyclic subgroups of $Q$,
which are not conjugate in $Q$, and
$$
\G = \langle \;  Q, y \; \|  \; y p_1 y^{-1} = p_2  \;  \rangle
$$
then $Q$ embeds in $\G / \G^n$.

Let $A \in \G$, $A$ not conjugate in $\G$ to an element of $Q$, and $A^{-1}
q_0 A = q_1$, $A^{-1} q_1 A = q_2$ in $\G$, where $q_0, q_1, q_2 \in Q$.
Without loss of generality, we can suppose that $A$ is a word of minimal
length (that is, $A = B$ in $\G$ implies $|A| \le |B|$). Then, using the
assumption that subgroups $\langle p_1 \rangle$, $\langle p_2 \rangle$ are
not conjugate in $Q$ and the remark that an equality of the form $q p^k
q^{-1} = p^\ell \neq 1$ in $Q$ implies $p^k = p^\ell$, we can obtain from
the equalities $A^{-1} q_0 A = q_1$, $A^{-1} q_1 A = q_2$ that $q_0 =q_1$.
Therefore, $q_0 \in \F(A)$, the automorphism $\kp_A$ is trivial and the
group $\la \kp_A, \F(A) \ra = \la \F(A) \ra$ has exponent $n$. By Theorem,
$Q$ naturally embeds in $\G / \G^n$ and Corollary 2 is proved. \qed
\smallskip

{\em Acknowledgements.}  The author  is grateful  to the referee
for finding a mistake in the proof of Corollary 1 and bringing
K.V. Mikhajlovskii's paper \cite{M94} to author's attention. The
author also wishes to thank A.Yu. Ol'shanskii for his kind
comments on results of this article and for pointing out that some
of them could also be proved using the techniques developed in his
joint with M.V. Sapir paper \cite{OS02} (see \cite{OS02} for
details).


\begin{thebibliography}{[HNN49]}


\bibitem[HNN49]{HNN49}  G. Higman, B.H. Neumann and H. Neumann,
{\em Embedding theorems for groups}, J. London Math. Soc.
{\bf 24}(1949), 247--254.

\bibitem[I94]{I94}
S.V. Ivanov,  {\em  The free Burnside groups of sufficiently
large exponents},  Internat. J. Algebra and Comp.  {\bf 4}(1994), 1--308.

\bibitem[I00]{I00}
S.V. Ivanov,  {\em  On finitely presented groups given
by periodic relators}, J. Group Theory  {\bf 3}(2000), 95--99.

\bibitem[I02]{I02}
S.V. Ivanov, {\em Weakly finitely presented infinite
periodic groups}, Contemporary Math.  {\bf 296}(2002), 139--154.

\bibitem[IO96]{IO96}
S.V. Ivanov and A.Yu. Ol'shanskii,
{\em Hyperbolic  groups and their  quotients of  bounded exponents},
Trans.  Amer. Math. Soc. {\bf 348}(1996), 2091--2138.

\bibitem[IO97]{IO97}
S.V. Ivanov and A.Yu. Ol'shanskii,
{\em On finite and locally finite subgroups of free Burnside
groups of large even exponents},  J.  Algebra  {\bf  195}(1997), 241--284.

\bibitem[LS77]{LS77}
R.C. Lyndon and P.E. Schupp, {\em Combinatorial group theory}, Springer-Verlag, 1977.

\bibitem[M94]{M94}
K.V. Mikhajlovskii, {\em Some generalizations of HNN-extensions in the
periodic case}, $\#$~1063-B94, VINITI, Moscow, 1994 (this is kept in Depot of
VINITI, Moscow, and is available upon request) 59pp. (in Russian)

\bibitem[O82]{O82a}
A.Yu. Ol'shanskii, {\em  On the Novikov-Adian theorem}, Mat. Sbornik {\bf
118}(1982), 203--235.

\bibitem[O89]{O89}
A.Yu. Ol'shanskii,  {\em Geometry of Defining Relations in Groups}, Nauka, Moscow, 1989;
English translation in Math. and Its Applications (Soviet series) {\bf 70} (Kluwer Acad.
Publishers, 1991).

\bibitem[OS02]{OS02} A.Yu. Ol'shanskii and M.V. Sapir,
{\em Non-amenable finitely presented torsion-by-cyclic groups}, Publ. IHES
{\bf  96}(2002), to appear.

\end{thebibliography}
\end{document}